\documentclass{amsart}
\usepackage{amsmath}
\usepackage{graphicx}
\usepackage{amsfonts}
\usepackage{amsthm}
\usepackage{amssymb}
\theoremstyle{plain}
\newtheorem{theorem}{Theorem}

\newtheorem{proposition}{Proposition}

\theoremstyle{definition}

\newtheorem{definition}{Definition}
\newtheorem{assumption}{Assumption}
\newtheorem{example}{Example}

\newtheorem{remark}{Remark}

\theoremstyle{remark}

\newcommand{\ds}{\displaystyle}

\newcommand{\dom}{\operatorname{dom}}
\newcommand{\Tan}{\operatorname{T}}
\newcommand{\bl}{\mathbf{l}}

\numberwithin{equation}{section}
\begin{document}
\title[Absolute continuity of Jump-Diffusions]{Absolutely continuous laws of Jump-Diffusions in finite and infinite dimensions with applications to mathematical Finance}
\author{Barbara Forster, Eva L\"utkebohmert, Josef Teichmann}
\address{Vienna University of Technology, Department of mathematical methods in
Economics, Research Group e105 Financial and Actuarial Mathematics,
Wiedner Hauptstrasse 8-10, A-1040 Wien, Austria}
\email{bforster@fam.tuwien.ac.at, eva.luetkebohmert@uni-bonn.de}
\email{jteichma@fam.tuwien.ac.at}
\thanks{The first author acknowledges the support from the FWF-project P15889 'Utility Maximization in Incomplete Financial Markets' and from the FWF-Wissenschaftskolleg W 8 'Differential Equation Models in Science and
Engineering'. The third author acknowledges the support from the RTN network
HPRN-CT-2002-00281 and from the FWF-grants Y328 and Z36. All authors are grateful to two unknown referees who helped to considerably improve the presentation and contents of this paper. In particular we are grateful for pointing out Corollary 3.25 of \cite{kusstr:85}}

\begin{abstract}
In mathematical Finance calculating the Greeks by Malliavin \newline weights has proved to be a numerically satisfactory procedure for finite-dimen\-sional It\^{o}-diffusions. The existence of Malliavin weights relies on absolute continuity of laws of the projected diffusion process and a sufficiently regular density. In this article we first prove results on absolute continuity for laws of projected jump-diffusion processes in finite and infinite dimensions, and a general result on the existence of Malliavin weights in finite dimension. In both cases we assume H\"ormander conditions and hypotheses on the invertibility of the so-called linkage operators. The purpose of this article is to show that for the construction of numerical procedures for the calculation of the Greeks in fairly general jump-diffusion cases one can proceed as in a pure diffusion case. We also show how the given results apply to infinite dimensional questions in mathematical Finance. There we start from the Vasi\v{c}ek model, and add -- by pertaining no arbitrage -- a jump diffusion component. We prove that we can obtain in this case an interest rate model, where the law of any projection is absolutely continuous with respect to Lebesgue measure on $\mathbb{R}^M $.
\end{abstract}
 \maketitle

\section{Introduction}

We shall consider in this article the question whether the law of
$\bl(X_t^x)$, for a finite dimensional projection $ \bl: H \to
\mathbb{R}^M,$ is absolutely continuous with respect to Lebesgue
measure on $\mathbb{R}^M,$ where $X_t^x$ is the solution of the
stochastic (partial) differential equation
\begin{align}
dX_{t}^{x}  &  =(AX^x_{t^-} + \alpha(X_{t^-}^{x}))dt+\sum_{i=1}^{d}V_{i}(X_{t^-}^{x})dB_{t}^{i}
+\sum_{j=1}^{m}\delta_{j}(X_{t^-}^{x})dL_{t}^{j},\\
X_{0}^{x}  &  =x \in H,
\end{align}
and $H$ is a possibly infinite dimensional separable Hilbert space.
We refer to the previous equation loosely speaking as a
jump-diffusion on the Hilbert space $H$ pointing out that the
involved L\'evy processes are of finite type. For sake of simplicity we shall always work with the c\`agl\`ad integrand $ t \mapsto X^x_{t^-} $, even though for the $ dt $ and $ dB_t $ integrals this is superfluous. In the infinite dimensional setting we are not aware of results on absolute continuity of the projected process in the jump diffusion case. Related work has been done in \cite{BelDal94} for the construction of first variation processes. In the diffusion case we refer to the work \cite{BauTei05} and the references therein, and in particular to the recently published inspiring results of Jonathan Mattingly, see for instance \cite{BakMat07}. We point out that we deal here with stochastic partial differential equations (SPDEs) and stochastic differential equations (SDEs) at the same time, where the latter case appears in this setting when the state (Hilbert) space $H$ is finite dimensional.

Very satisfying results in the finite dimensional setting with L\'evy processes
of infinite type have been obtained by in \cite{Cass07} through a
generalization of the Norris Lemma to L\'evy processes. These results have been built upon our results presented in this work for the finite activity
case (see Section 7). Substantial work with respect to absolute continuity and smoothness of the density has been already published in the eighthies, where the most prominent ones are \cite{Bis83} and \cite{Bich87}. Therein several questions of extension of hypo-ellipticity results (and Malliavin Calculus) to jump-processes are discussed and completely solved, however, the problem of a hypo-elliptic diffusion together with a finite-activity jump-structure remained open. It has to be pointed out here that -- in contrast to \cite{Bich87} -- we do not need any extension of Malliavin Calculus to jump processes for our results (see also the discussion in Remark \ref{product-structure}). This gap was filled by the announced results of \cite{Zhou92}, but several proofs therein are extremely short. Recently -- motivated by questions from financial mathematics, see Section 8 for an outline of the problem -- there has been increasing interest in those results, see the the works \cite{Bal07}, \cite{DavJoh06} and \cite{El-Kha04} and the references therein. This article aims to work out the most general finite activity case under H\"ormander conditions on the diffusion part.

From the point of view of existence and uniqueness for jump-diffusions in infinite dimensions our main reference is \cite{FilTap07}, and the references therein. Since we consider jump-diffusions as concatenated diffusions on Hilbert space we mention \cite{DaPZab92} as the main reference for existence and uniqueness results, but also \cite{BelDal90} and \cite{CarRoz99} for many interesting constructions and ideas.

There are two applications added to this work. The first one is the HJM-equation (as presented in \cite{FilTap07}), where we show that the innocent Vasi\v{c}ek model (see for instance the seminal work of \cite{BjoSve01}) with a certain jump structure triggered by a one-dimensional Poisson process yields -- under no-arbitrage assumptions -- a model, where not only no finite dimensional realizations do exist, but where every projection into a finite dimensional subspace admits a density (compare with the notion of generic interest rate evolutions from \cite{BauTei05}). The second application is concerned with concrete formulas for the calculation of Malliavin weights. There our message is that one can think Poisson-trajectory-wise, i.e. the results from \cite{gobmun:05} or \cite{Fou99} can be literally applied by replacing the diffusion process by the respective jump-diffusion process.

When we analyse jump-diffusions with values in Hilbert spaces then loosely speaking the following facts hold true:

\begin{itemize}
\item between two consecutive jumps of the jump-diffusion we are given an
ordinary diffusion.

\item at a jump we add to the left limit the jump size (which
usually depends on the left limit, too). In \cite{Pro90},
Chapter V.10, Hypothesis (H3), this operation is formalized by the so
called \emph{linkage operators} $ x \mapsto x + \mu \delta^j(x) $, which encode what happens at a jump of size $ \mu $ at $ x $. We shall apply this notion
here, too.
\end{itemize}

Hence the following picture arises:

\begin{itemize}
\item In order to obtain absolute continuity of the projected diffusion process, we need the H\"ormander condition to be in force. Otherwise we cannot expect -- conditioned on the event of positive probability that no jump occurs -- that the law of $ \bl(X_t^x) $ is absolutely continuous with respect to Lebesgue measure. 

\item In order to preserve absolute continuity we need the linkage operators to be invertible in a proper sense.
\end{itemize}
\begin{remark}
Both conditions are 'sine qua non', since it is easy to imagine counterexamples.
\end{remark}
In Section 2 we fix the general setting of this article. We shall deal with L\'evy processes of
finite type as drivers of the stochastic differential equations,
even though we believe that one should be able to prove similar
results in the case of many small jumps, too. We also state the main
assumptions of this work in Section 2 for later use.

In Section 3 we prove a ``folklore'' decomposition theorem, which tells that solving a jump-diffusion S(P)DE
is the same as solving associated diffusion S(P)DEs and concatenating
the solutions by linkage operators at jumps. In Section 4 we show that we can
also prove results on first variation processes in the spirit of the
decomposition theorem. We prove that under our analytic
requirements there is in fact a sufficiently regular first variation
process. 

In Section 5 we show by means of Malliavin calculus for a
$d$-dimensional Brownian motion that the law of a projected
jump-diffusion is absolutely continuous with respect to Lebesgue
measure. In Section 6 we introduce a class of examples from mathematical Finance,
where we see very directly the phenomenon of absolute continuity arising from the introduction of jumps and the no-arbitrage condition. This shows again that finite dimensional realizations, as constructed for instance in \cite{bjogom:99}, are a rare case in infinite dimensions. In Section 7 we restrict to the finite dimensional setting to show that
the density of the absolutely continuous law is in fact smooth by proving that the inverse of the covariance matrix has $p$-th moments for all $ p \geq 1$. In Section 8 we apply the invertibility of the covariance matrix to the calculation of Greeks. Section 9, an Appendix, shows an important estimate implicitly present in the Norris Lemma as presented in D.~Nualart's book \cite{Nua95}. A similar result (which could be directly used in Section 7 for the proof of the main theorem) can be found in \cite[Corollary 3.25]{kusstr:85}. The article \cite{kusstr:85} is most likely the source of the first appearance of the precise polynomial time-dependence in the estimate of the $L^p$-norm of the inverse of the covariance matrix. We have been choosing here the path via the Norris lemma, we explain the estimate by re-doing its proof in the Appendix.

\section{Setting and Assumptions}

Let $(\Omega,\mathcal{F},P,(\mathcal{F}_{t})_{t\geq0})$ be a filtered
probability space where the filtration $(\mathcal{F}_{t})_{t\geq0}$ satisfies
the usual conditions. Let $(B_{t})_{t\geq0}$ be a $d$-dimensional Brownian
motion and $(L_{t}^{j})_{t\geq0}$, $j = 1, \ldots, m$ be $m$ independent
compound Poisson processes given by
$$
L_t^j := \sum_{k = 1}^{N_t^j} Z_k^j,\label{poisson-process}
$$
where $N_t^j$ denotes a Poisson process with  jump intensity $\widetilde{\lambda_{j}}>0$
and $Z^j = (Z_k^j)_{k \geq 1}$ is an i.i.d.~sequence of random variables with
distribution $\mu_{j}$ for $j=1,\dots,m$, such that each $\mu_{j}$ admits all
moments. The compensated compound Poisson
process reads as follows,
$$
L_t^j - E(L_t^j) = L_t^j - \lambda_j t,\label{compensator}
$$ where $\lambda_j =
E(Z_1^j) \widetilde{\lambda_j}$ is the average jump size times the jump rate.

We could equally take an $ \mathbb{R}^m $-valued L\'evy process of finite type, i.e.~introduce a dependence structure between the jumps of the components and all theorems would equally hold true with slightly modified proofs, but we believe that this generalisation does not bring further insight.

We assume that all sources of randomness are mutually independent
and that the filtration $(\mathcal{F}_{t})_{t\geq0}$ is
the natural filtration with respect to $(B_{t},L_{t}^{1},\dots,L_{t}
^{m})_{t\geq0}$.
Let $ H $ be a separable Hilbert space. We fix furthermore a
strongly continuous semi-group $ S $ on $H$ with generator $ A $.
Let $\alpha,V_{1},\dots,V_{d}$, the diffusion vector fields, and
$\delta _{1},\dots,\delta_{m}$, the jump vector fields, be
$C^{\infty}$-bounded on $H$, that is, the vector fields are
infinitely often differentiable with bounded partial derivatives of
all proper orders $ n \geq 1 $. We consider the mild c\`{a}dl\`{a}g
solution $(X_{t}^{x})_{t \geq 0}$ of a stochastic differential
equation
\begin{align}
dX_{t}^{x}  &  =(AX^x_{t^-} + \alpha(X_{t^-}^{x}))dt+\sum_{i=1}^{d}V_{i}(X_{t^-}^{x})dB_{t}^{i}
+\sum_{j=1}^{m}\delta_{j}(X_{t^-}^{x})dL_{t}^{j},\\
X_{0}^{x}  &  =x \in H.
\end{align}
See \cite{FilTap07} for all necessary details on existence and uniqueness of the previous equation.

The previous conditions are slightly more than standard for
existence and uniqueness of mild solutions, i.e. in \cite{FilTap07} the authors need Lipschitz conditions on the vector fields, whereas we assume them to be $C^{\infty}$-bounded. 
In order to speak about absolute continuity of projections to $\mathbb{R}^M$ we shall need more assumptions, in particular for conclusions drawn from the geometry of the given vector
fields $ \alpha, V_1,\dots,V_d,\delta_1,\ldots,\delta_m$ several
quite strong analytic requirements are necessary. We group the
assumptions in three groups and indicate in each section, which
assumptions we shall need.

Let $ \bl :H \to \mathbb{R}^M $ be a projection, then we want to
know whether the law of $ \bl(X_t^x) $ is absolutely continuous with
respect to Lebesgue measure and if -- in case -- the density is
smooth. Following the short discussion in the introduction we need the
H\"ormander conditions to be in force and we need to suppose
invertibility on linkage operators.

We apply the following notations for Hilbert spaces $ \dom(A^k) $,
$$
\begin{array}{ccl}
\dom(A^k)&:=&\ds \left\{h\in H| \, h\in \dom(A^{k-1})
\hspace{2mm}\mbox{and}\hspace{2mm}
A^{k-1}h\in \dom(A)\right\},\\
||h||^2_{\dom(A^k)}&:=&\ds \sum_{i=0}^k ||A^ih||^2,\\
\dom(A^{\infty})&=&\ds \bigcap_{k\geq 0} \dom(A^k),
\end{array}
$$
which we need in order to specify the analytic conditions.

\begin{assumption}\label{ass analytic conditions}
We assume that the generator $A$ of $S$ generates a \emph{strongly
continuous group}. We assume furthermore that
$\alpha,V_{1},\dots,V_{d}$, the diffusion vector fields, and $\delta
_{1},\dots,\delta_{m}$, the jump vector fields, are
$C^{\infty}$-bounded on the Hilbert spaces $\dom(A^k)$ for $ k \geq
0 $, that is, the vector fields are infinitely often differentiable
with bounded partial derivatives of all proper orders $ n \geq 1 $
on the Hilbert space $ \dom(A^k) $ for $ k \geq 0 $.
\end{assumption}

\begin{assumption}\label{ass hoermander condition}
We take Assumption \ref{ass analytic conditions} for granted, i.e.
we can consider all vector fields on the space $ \dom(A^k) $ for $
k=0,\ldots,\infty $. For a proper statement of the H\"ormander
condition we apply the ``geometrically relevant'' drift
$$
V_0(x) = Ax + \alpha(x) - \frac{1}{2} \sum_{i=1}^d \Tan V_i(x) \cdot V_i(x)
$$
for $ x \in \dom(A)$ and call $V_0$ the Stratonovich drift of the diffusion. Recall the (tangent) directional derivative operator $\operatorname{T}$ defined through
$$\operatorname{T}V(x) \cdot v = \frac{d}{d\epsilon}\Big |_{\epsilon = 0} V(x + \epsilon v).$$
Lie brackets can only be calculated on the Fr\'echet space $
\dom(A^{\infty}) $ and there we formulate the H\"ormander condition.
We assume that the distribution $\mathcal{D}(x)$ generated by the
vector fields
\begin{gather}\label{hoermander condition}
V_{1}(x),\dots, V_{d}(x), \quad [V_{i}(x) ,V_{j}(x)] \; (i, j
=0, 1,\dots,d),\\
[V_i(x),[V_j(x), V_k(x)]] \; (i, j, k
=0, 1, \ldots, d), \quad \ldots \ldots \nonumber
\end{gather}
is dense in $H$ for one $x\in \dom(A^{\infty})$.
\end{assumption}

\begin{assumption}\label{ass jump condition}
We assume that the inverse of $x\mapsto x+z \delta_{j}(x)$ exists
and is $C^{\infty}$-bounded on each $ \dom(A^k) $ for
$z\in\operatorname*{supp}(\mu_{j})$, $j=1,\dots,m$ and $ k \geq 0 $
(recall that $\mu_j$ was the distribution of the random variable
$Z_j$).
\end{assumption}

\begin{remark}
As far as Assumption \ref{ass analytic conditions} is concerned we
do believe that the assertions of this paper also hold true for (most) strongly continuous
semigroups. A proof based on an application of the Szek\H{o}falvi-Nagy theorem can be found in the recent preprint \cite{tei08}, therefore we could replace the assumption that $ A $ generates a strongly continuous group by the assumption that $ A $ generates a pseudo-contractive strongly continuous semigroup. This includes most of the second order partial differential operators. However, for this paper we do always assume the group property for the sake of simplicity.
\end{remark}
\begin{remark}
The H\"ormander condition could not be formulated without the analytic part of Assumption \ref{ass analytic conditions}.
\end{remark}

\begin{example}\label{banach-map-vector-fields}
In order to show examples of vector fields, which are $ C^{\infty}$-bounded on $ \dom(A^k) $ consider the following structure. Let $H$ be a separable Hilbert space and $ A $ the generator of a strongly continuous semigroup. We know that $ \dom(A^{\infty}) $ is a Fr\'echet space and an injective limit of the Hilbert spaces $\dom(A^k) $ for $ k \geq 0$. Following the analysis as developed in \cite{Fil03} (see also \cite{KriMic97} and \cite{Ham82} were the analytic concepts have been originally developed), we can consider the following vector fields $ V: U \subset H \to \dom(A^{\infty}) $. If $V$ is smooth in the sense explained in \cite{Fil03} and has the property that its derivatives of proper order $ n \geq 1 $ are bounded on $ U \subset H $, then $ V$ is obviously a $ C^{\infty}$-bounded vector field and additionally $ V|_{\dom(A^{\infty})} $ is a Banach-map-vector field in the sense of \cite{Fil03}. Such vector fields constitute a class, where Assumptions 1--3 can be readily checked.
\end{example}

\section{Decomposition Theorem for Jump-Diffusions on Hilbert Spaces}

In order to properly understand how to apply the Malliavin calculus,
we state the following rather obvious structure theorem on
jump-diffusions, which simply takes into account that stochastic
integration with respect to the Poisson process follows the rules of
Lebesgue-Stieltjes integration (see for instance \cite{Pro90} for a
general exposition). Here we only need that the vector fields are $
C^{\infty}$-bounded on $H$ in order to guarantee existence and
uniqueness of the respective equations.

\begin{theorem}
Let $(\Omega,\mathcal{F},P,(\mathcal{F}_{t})_{t\geq0})$ be a
filtered probability space and let $(B_{t})_{t\geq0}$ be a
$d$-dimensional Brownian motion and $(L_{t}^{j})_{t\geq0}$ be $m$
independent compound Poisson processes for $j=1,\dots,m$, such that the
filtration is the natural filtration with respect to
$(B_{t},L_{t}^{1} ,\dots,L_{t}^{m})_{t\geq0}$. Let $ S $ be a
strongly continuous semigroup with generator $ A $ on $H$. Let
$\alpha,V_{1},\dots,V_{d},$ the diffusion vector fields, and
$\delta_{1},\dots,\delta_{m},$ the jump vector-fields, be
$C^{\infty}$-bounded on $ H $ and consider the c\`{a}dl\`{a}g
solution $(X_{t}^{x})_{0\leq t\leq T}$ of a stochastic differential
equation
\begin{align}
dX_{t}^{x}  &  =(AX_{t^-}^{x} + \alpha(X_{t^-}^x))dt+\sum_{i=1}^{d}V_{i}(X_{t^-}^{x})dB_{t}^{i}
+\sum_{j=1}^{m}\delta_{j}(X_{t^-}^{x})dL_{t}^{j},\label{equ1}\\
X_{0}^{x}  &  =x. \label{equ2}
\end{align}
Let $\eta$ denote a piecewise constant, c\`{a}dl\`{a}g trajectory
$\eta: \mathbb{R}_{\geq 0} \to \mathbb{R}^{m}$ of the compound
Poisson process $L$ with finitely many jumps on compact intervals
and starting at $0$. We consider
\begin{align}
dY_{s,t}^{x,\eta}  &  =\Big(A Y_{s,t^-}^{x,\eta} + \alpha(Y_{s,t^-}^{x,\eta}) \Big)dt+\sum_{i=1}^{d}V_{i}(Y_{s,t^-}^{x,\eta}
)dB_{t}^{i}\label{equ1 along eta}\\
&  +\sum_{j=1}^{m}\delta_{j}(Y_{s,t^-}^{x,\eta})d\eta^{j}(t),\nonumber\\
Y_{s,s}^{x,\eta}  &  =x. \label{initial value}
\end{align}
Then $(Y_{s,t}^{x,\eta})_{s \geq 0} $ can be given explicitly in
terms of the jump times $ \tau_n $ of $\eta$ for $ n \geq 0 $ and the
diffusion process between two consecutive jumps:
\begin{align*}
Y_{t}^{x,\eta}  &  :=Y_{0,t}^{x,\eta}\text{ for }0\leq t<\tau_{1}\\
Y_{t}^{x,\eta}  &
:=Y_{\tau_{1},t}^{y,\eta}\big|_{y=Y_{0,\tau_{1}^-}^{x,\eta}+\sum
_{j=1}^{m}\delta_{j}(Y_{0,\tau_{1}^-}^{x,\eta})\Delta\eta^{j}(\tau_{1})\text{ }
}\text{ for }\tau_{1}\leq t<\tau_{2}\\
&  \vdots &\\
Y_{t}^{x,\eta}  &  :=Y_{\tau_{n-1},t}^{y,\eta}\big|_{y=Y_{0,\tau_{n-1}^-}^{x,\eta}
+\sum_{j=1}^{m}\delta_{j}(Y_{0,\tau_{n-1}^-}^{x,\eta})\Delta\eta^{j}(\tau_{n-1})\text{
}}\text{ for }\tau_{n-1}\leq t<\tau_{n}.
\end{align*}
Here we write $\Delta\eta(t):=\eta(t)-\eta(t^-)$ for
$t \geq 0$. We define the process $(Y_{t}^{x,L})_{t\geq 0}$ by inserting the compound Poisson process
$L$ for $\eta$ in $(Y_{t}^{x,\eta })_{t \geq 0}$. The resulting process $(Y_{t}^{x,L})_{t \geq 0}$
is then indistinguishable from $(X_{t}^{x})_{t \geq 0}$.
\label{structure theorem}
\end{theorem}

\begin{proof}
For the proof we refer to \cite{Pro90}, Chapter V.10, Theorem
57, in particular with respect to the conditioning on the jump
part. The proof remains unchanged in the infinite dimensional setting, see \cite{FilTap07} for the existence and uniqueness proof on separable Hilbert spaces.
\end{proof}

\begin{remark}\label{remark_on_eta}
For future use we shall always assume that the first jumping time of
$\eta$ is strictly positive $ \tau_1 > 0$ and that each time
corresponds to the jump of exactly one coordinate process $L^j$,
which is both true for almost all trajectories of the compound Poisson
process $ L $. Notice that the dependence of
$(Y_{t}^{x,\eta})_{t\geq0}$ on the jump times of $\eta$ is
continuous, but certainly not smooth since the jump times are inserted instead of time of a hypo-elliptic diffusion process.
\end{remark}

\begin{remark}\label{product-structure}
Notice that one can also interpret the result in the following way:
consider the solution $(X_{t}^{x})_{t\geq0}$ of equation
(\ref{equ1}) as an element of $L^{2}(\Omega_{1}\times\Omega_{2};H)$,
where $\Omega_{1}$ carries the Brownian motion part (with natural
filtration), $\Omega_{2}$ carries the Poisson part (with natural
filtration) and $\Omega_1 \times \Omega_2$ is equipped with the
respective product $\sigma$-algebra. Then we know by Fubini's
Theorem that
\[
L^{2}(\Omega_{1}\times\Omega_{2};H)=L^{2}(\Omega_{2};L^{2}
(\Omega_{1};H)).
\]
The previous theorem only clarifies the jump-diffusion structure
of the dependence on $\Omega_{2}$. In other words, between jumps
we have ordinary diffusions, and at a jump we link by linkage
operators.
\end{remark}

\section{First Variation Processes}

In order to calculate Malliavin derivatives, which is crucial for arguments on absolute continuity, we need precise statements on first variation processes of jump diffusions. For later purposes, but also in order to see results on the inverse of the first variation process easily, we write our equations in the Stratonovich notation. This is not innocent in infinite dimensions, since mild solutions are in general \textbf{not} semi-martingales and therefore the Stratonovich notation fails to be applicable in general. However, by Assumption \ref{ass analytic conditions} we are able to determine whether we are given a semi-martingale or not by analysing the initial value of the process. Indeed -- for fixed $ k \geq 0 $, if $ x \in \dom(A^{k+1}) $, then there is a mild solution taking values in $ \dom(A^{k+1}) $ of the It\^o stochastic differential equation. However, this solution process has to coincide -- by uniqueness -- with the solution process obtained by considering the same equation on $ \dom(A^k) $ with an initial value in $ \dom(A^{k+1}) $. Therefore the mild solution in $ \dom(A^{k+1}) $ is a strong solution in $ \dom(A^k) $. Therefore we assume Assumptions \ref{ass analytic conditions} and \ref{ass jump condition} to be in force in this section.

By the previous arguments for the given stochastic differential equation (\ref{equ1}), we can switch to Stratonovich
notation for $ x \in \dom(A) $, and obtain
$$
dX_{t}^{x} =V_0(X_{t^-}^{x})dt+\sum_{i=1}^{d}V_{i}(X_{t^-}^{x}) \circ dB_{t}^{i}
+\sum_{j=1}^{m}\delta_{j}(X_{t^-}^{x})dL_{t}^{j}
$$
with the Stratonovich drift given by
\[
V_{0}(x):=Ax+\alpha(x)-\frac{1}{2}\sum_{i=1}^{d}\Tan V_{i}(x)\cdot V_{i}(x)
\]
for $ x \in \dom(A) $. Recall the tangent (derivative) operator $\operatorname{T}$
$$\operatorname{T}V(x) \cdot v = \frac{d}{d\epsilon}\Big |_{\epsilon = 0} V(x + \epsilon v).$$

We do also consider the stochastic differential equation (\ref{equ1 along eta}) with respect to one trajectory $\eta$ and switch to Stratonovich notation, too. The following theorem states the result on the first variation process along one trajectory $ \eta $, which yields in the sequel the same result by inserting the compound Poisson process $L$ for $ \eta $. Notice that the trajectory $ \eta $ is such that the first jumping time is strictly positive and that at each jumping time $ \tau_n $ for $ n \geq 1 $ only one coordinate jumps.

\begin{theorem}
Assume Assumptions \ref{ass analytic conditions} and \ref{ass jump condition}. We fix $ k \geq 0 $.
The \emph{first variation process} ${(J_{s \to t}(x,\eta))}_{t \geq s}$ associated with
${(Y_t^{x, \eta})}_{t \geq 0} $ on $\dom(A^k)$ is well defined and satisfies the
stochastic differential equation
\begin{align}
dJ_{s\to t}(x,\eta) \cdot h & =\Big(A J_{s\to
t-}(x,\eta) \cdot h+ \Tan \alpha(Y_{s,t^-}^{x,\eta})\cdot
J_{s\to t^-}(x,\eta)\cdot h\label{firstvariation} \Big)dt\nonumber\\
&  +\sum_{i=1}^{d} \Big(\Tan V_{i}(Y_{s,t^-}^{x,\eta})\cdot J_{s\to t^-}
(x,\eta)\cdot h \Big) dB_{t}^{i}\nonumber \\
&  +\sum_{j=1}^{m} \Big(\Tan \delta_{j}(Y_{s,t^-}^{x,\eta})\cdot
J_{s\to t^-}(x,\eta)\cdot h \Big) d\eta^{j}(t),\nonumber\\
J_{s\to s}(x,\eta) \cdot h&  =h,
\end{align}
for $h,x\in \dom(A^k)$ and $t\geq s.$ The It\^{o} equation has a unique global mild solution for $h,x\in\dom(A^k)$ and $J_{s\to t}(x,\eta)$ defines a continuous linear
operator on $\dom(A^k)$, which is invertible if
$x\in\dom(A^{k+1})$. \\
The Stratonovich equation on $\dom(A^k)$ in turn is only
well-defined for $h,x\in\dom(A^{k+1})$. We apply the (formal) notation here
$$
\Tan V_0(x) v = Av+\Tan \alpha(x)v-\frac{1}{2} \Tan (x \mapsto \sum_{i=1}^{d}\Tan V_{i}(x)\cdot V_{i}(x))v
$$
for $ x \in \dom(A) $ and $ v \in \dom(A) $.
\begin{align}
dJ_{s\to t}(x,\eta) \cdot h &
= \Big( \Tan V_{0}(Y_{s,t^-}^{x,\eta})\cdot
J_{s\to t^-}(x,\eta)\cdot h \Big)dt\label{firstvariationStratonovich}\nonumber \\
& \hspace{2mm} +\sum_{i=1}^{d}\Big( \Tan V_{i}(Y_{s,t^-}^{x,\eta})\cdot J_{s\to t^-}
(x,\eta)\cdot h \Big) \circ dB_{t}^{i} \\
&  \hspace{2mm} +\sum_{j=1}^{m}\Big( \Tan \delta_{j}(Y_{s,t^-}^{x,\eta})\cdot
J_{s\to t^-}(x,\eta)\cdot h \Big) d\eta^{j}(t),\nonumber\\
J_{s\to s}(x,\eta) \cdot h&  =h,\nonumber
\end{align}
for $h,x\in\dom(A^{k+1})$ and $t\geq s.$ The adjoint of the
inverse
$$ Z_t^{x,h} := \left(J_{s\to t}(x,\eta)^{-1}\right)^{*} \cdot h,$$
if it exists, should satisfy the following Stratonovich equation at
the point $ x $ in direction $ h $,
\begin{align}
dZ_t^{x,h}& =-\Big(
\Tan V_0(Y_{s,t^-}^{x,\eta})^{*}\cdot Z_{t^-}^{x,h}\Big)dt -\sum_{i=1}^{d} \Big( \Tan V_{i}(Y_{s,t^-}^{x,\eta })^{*} \cdot Z_{t^-}^{x,h} \Big) \circ dB_{t}^{i} - \nonumber\\
&\hspace{2mm}-\sum_{j=1}^{m}\Big( \Tan \delta_{j}(Y_{s,t^-}^{x,\eta})^{*} \cdot Z_{t^-}^{x,h} \Big) d\eta^{j}(t) \label{firstvariationinv}\\
&\hspace{2mm}+\sum_{j=1}^{m}\Big(\left( \Tan \delta_{j}(Y_{s,t^-}^{x,\eta})^{2}\right)^{*}\cdot \left(\left(  id_{H}
+\Delta\eta^{j}(t) \Tan \delta_{j}(Y_{s,t^-}^{x,\eta})\right)
^{-1}\right)^{*}\cdot \nonumber\\
&\hspace{1cm} \cdot Z_{t^-}^{x,h} (\Delta\eta^{j}(t))^{2}\Big),\nonumber
\end{align}
for $h,x\in\dom(A^{k+1})$ and $t\geq s\geq 0$ (Here we
applied the notions of \cite{Pro90}).
\end{theorem}
\begin{remark}
The completely analogous theorem holds when we replace $ \eta $ by a
compound Poisson process $L$. We do not state this theorem again,
but we point out that we even have moment estimates for the
respective processes, which is the only additional relevant
information. To be precise, the first variation process $ J_{s \to
t}(x)\cdot h $, which equals $ J_{s \to t}(x,L) $ by construction,
has bounded second moments by \cite{FilTap07}.
\end{remark}
\begin{proof}
Under our Assumption \ref{ass analytic conditions} the regularity in
the initial values is clear by well-known results from
\cite{DaPZab92} and the chain rule on Hilbert spaces (recall that
the linkage operators are smooth).  We are allowed to pass to the
Stratonovich decomposition since we integrate semi-martingales by
It\^{o}'s formula on Hilbert spaces for $x,h\in \dom(A^{k+1})$ due
to the arguments of \cite{BauTei05}: the core assertion is here that
we can replace $ H $ by each $ \dom(A^k) $ for some $ k \geq 0 $,
which means in turn if we start in $ \dom(A^{k+1})$ and obtain a
mild solution there, it is indeed a strong solution considered on $
\dom(A^k) $, for $ k \geq 0 $. It remains to show the invertibility
results on the respective first variation processes.

Left invertibility of the first variation $J_{s\to t}(y,.)$ follows
by It\^{o}'s formula since we have c\`{a}dl\`{a}g trajectories with
finitely many jumps. Calculating the semi-martingale decomposition
of ${(Z_t^{y,.})}^{\ast} \cdot J_{0\to t}(x,\eta)$ given by
equations (\ref{firstvariationStratonovich}) and
(\ref{firstvariationinv}) yields the result
$$
{(Z_t^{x,.})}^{\ast} J_{0\to t}(x,\eta)=id_{\dom(A^k)}.
$$
Thus, the solution of equation (\ref{firstvariationinv}) is the left
inverse of $J_{s\to t}.$\\
We prove that the left inverse is also the right inverse by the same
reasoning as in the proof of Proposition 2 in \cite{BauTei05}.
Therefore, we choose an orthonormal basis $(g_i)_{i\geq 1}$ of
$\dom(A^k)$ which lies in $\dom(A^{k+1}).$ Then we can compute the
semi-martingale decomposition of
$$
\begin{array}{ccl}
\ds\sum_{i=1}^N \langle {(Z_t^{x,h_1})}^{\ast},g_i\rangle_{\dom(A^k)}\langle g_i, J_{s\to t}
(x,\eta)^{*}\cdot h_2\rangle_{\dom(A^k)}&=&\\[4mm]
\ds\sum_{i=1}^N\langle h_1,Z_t^{x,g_i}\rangle_{\dom(A^k)}\langle
J_{s\to t} (x,\eta)\cdot g_i,h_2\rangle_{\dom(A^k)},
\end{array}
$$
for $h_1,h_2\in \dom(A^{k+1})$ and $N\geq 1.$ Applying the
Stratonovich decomposition and by adjoining we can free the
$g_i\!$'s and pass to the limit, which yields vanishing finite
variation and martingale part. Hence
$$
\begin{array}{c}
\ds\langle J_{s\to t} (x,\eta) {(Z_t^{x,h_1})}^{\ast}, h_2\rangle_{\dom(A^k)}=\\[4mm]
\ds\lim_{N\to \infty}\sum_{i=1}^N\langle {(Z_t^{x,h_1})}^{\ast},g_i\rangle_{\dom(A^k)}\langle
g_i,J_{s\to t} (x,\eta)^{*}\cdot
h_2\rangle_{\dom(A^k)}\\[4mm]
=\langle h_1,h_2,\rangle_{\dom(A^k)},
\end{array}
$$
which is what a right inverse should satisfy.
\end{proof}

\section{Absolutely continuous laws in finite and infinite dimensions}

In this section we assume Assumptions \ref{ass analytic conditions},
\ref{ass hoermander condition} and \ref{ass jump condition}. We want
to determine by means of Malliavin calculus whether the law of $
\bl(Y_{t}^{x,\eta})$ is absolutely continuous with respect to
Lebesgue measure for $t>0$.

For details on Malliavin calculus see \cite{Mal97} and \cite{Nua95},
where in particular the derivative operator and the Skorohod
integral for Malliavin calculus with respect to a $d$-dimensional
Brownian motion are defined. Notice that we do not need a Malliavin calculus with respect
to the Poissonian trajectories, since we calculate Poisson-trajectory-wise.

Our first task is the calculation of the Malliavin derivative for a
fixed c\`{a}dl\`{a}g path $\eta$. In a second step, we consider the
composed problem, where we replace $ \eta $ by a compound Poisson
process $L$ as outlined before. Therefore, we first fix a piecewise
constant, c\`{a}dl\`{a}g trajectory $\eta: \mathbb{R}_{\geq 0} \to
\mathbb{R}^m$ of the process $(L^1_t,\ldots,L^m_t)_{t\geq 0}$ with
finitely many jumps on compact intervals starting at $0$.

\begin{theorem}\label{theorem2}
We take Assumptions \ref{ass analytic conditions}, \ref{ass
hoermander condition} and \ref{ass jump condition} for granted,
where $ x \in \dom(A^{\infty})$ denotes the point where the
H\"ormander condition (\ref{hoermander condition}) holds true. Let
${(Y_{t}^{x})}_{t\geq0}$ denote the unique c\`{a}dl\`{a}g solution
of equation (\ref{equ1 along eta}). Then for projections $\bl :H\to
\mathbb{R}^M$ the law of $\bl(Y_{t}^{x})$ is absolutely continuous
with respect to Lebesgue measure on $\mathbb{R}^M$ for $t>0$.
\end{theorem}

\begin{proof}
Fix $t>0.$ We are able to write the Malliavin derivative of $Y^x_t$ for each Poissonian trajectory $ \eta $,
\[
D_{s}^{i}(\bl\circ Y_{t}^{x,\eta})=\bl\circ J_{0\to t}(x)J_{0\to
s}(x)^{-1} V_{i}(Y_{s^-}^{x,\eta})1_{[0,t]}(s).
\]
We can calculate the \emph{Malliavin covariance matrix} $\gamma$ as
$$\langle \gamma(\bl\circ Y_t^{x, \eta}) \xi, \xi \rangle := \sum_{i = 1}^d
\int_0^t \langle \bl\circ J_{0\to t}(x)J_{0\to s}(x)^{-1}
V_{i}(Y_{s^-}^{x,\eta}), \xi \rangle^2 ds.$$ Consequently, the
covariance matrix $\gamma(\bl\circ Y_t^{x, \eta})$ can be calculated
in the usual way via the reduced covariance matrix
\[
\left\langle C_{t}\xi,\xi\right\rangle :=\sum_{i=1}^{d}\int_{0}^{t}
\left\langle J_{0\to s}(x)^{-1} V_{i}(Y_{s^-}^{x,\eta}),\xi
\right\rangle^{2}ds.
\]
through the relation 
$$\gamma(\bl\circ Y_t^{x, \eta}) =(\bl\circ J_{0 \to t}(x)) C_t {(\bl\circ J_{0 \to t}(x))}^{\ast}),
$$ 
where $ \ast $ denotes the adjoint operator with respect to the Hilbert space structures on $H$ and $ \mathbb{R}^M$. We assume $ n $ jumps of $ \eta $ on $ [0,t] $ and we denote by $0=\tau_0<\tau_1<\ldots<\tau_n\leq t$ the sequence of jump times of $\eta.$ For convenience, we denote the last point in
time $t$ by $\tau_{n+1}$, even if $ \tau_{n+1} = \tau_n $, which can in principle happen.
Hence we can decompose,
\[
\left\langle C_{t}\xi,\xi\right\rangle :=\sum_{k=0}^n\sum_{i=1}^{d}\int_{\tau_k}^{\tau_{k+1}}
\left\langle J_{0\to s}(x)^{-1} V_{i}(Y_{s^-}^{x,\eta}),\xi
\right\rangle^{2}ds=\sum_{k=0}^n \langle C_t^k\xi,\xi\rangle.
\]
Each of the summands determines a symmetric matrix $C_t^{k}$ and can
be interpreted as a reduced covariance matrix coming from a
diffusion between $\tau_k$ and $\tau_{k+1}$ with initial value $
Y^x_{\tau_k^-} $ for $k=0,\ldots,n$. We do not know whether the
H\"{o}rmander condition is true everywhere. Therefore, we do not
know whether $C_t^k$ is a positive definite operator for all
$k\geq0.$ From \cite{BauTei05}, Theorem 1, we do know, however, that $C_t^{0}$ is a positive definite operator and there exist null
sets $N_{0}$ such that on $N_{0}^c$ the matrix $C_t^{0}$ is
invertible. Hence the law of $(\bl\circ Y_t^{x,\eta})$ is absolutely
continuous with respect to the Lebesgue measure on $\mathbb{R}^M$,
since $J_{0\to t}(x)$ is invertible and therefore $\gamma(\bl\circ
Y_t^{x, \eta})$ has empty kernel (Theorem 2.1.2 in \cite{Nua95}, p.
86).
\end{proof}

\begin{remark}
The same conclusions hold for $Y_{t^-}^{x,\eta}$: notice that
$Y_{t}^{x,\eta }=Y_{t^-}^{x,\eta}$, if there is no jump at $t$.
Otherwise $Y_{t}^{x,\eta
}=Y_{t^-}^{x,\eta}+\sum_{j=1}^{m}\delta_{j}(Y_{0,t^-}^{x,\eta})\Delta\eta
^{j}(t)$, but invertible diffeomorphisms transform absolutely continuous laws in absolutely continuous ones.\end{remark}

Now we extend this theorem to the jump-diffusion process
$(X_t^x)_{t\geq 0}$, which is easy since -- conditioned on one
trajectory $ \eta $ -- we do have an absolutely continuous law and this
property is not perturbed by integration due to Fubini's Theorem.

\begin{theorem}
We take Assumptions \ref{ass analytic conditions}, \ref{ass
hoermander condition} and \ref{ass jump condition} for granted,
where $ x \in \dom(A^{\infty})$ denotes the point where the
H\"ormander condition (\ref{hoermander condition}) holds true. Let
${(X_{t}^{x})}_{t\geq0}$ denote the unique c\`{a}dl\`{a}g solution
of equation (\ref{equ1}). Then for projections
$\bl=(l_1,\ldots,l_k):H\to \mathbb{R}^M$ the law of $ \bl(X_{t}^x)$
is absolutely continuous with respect to Lebesgue measure on
$\mathbb{R}^M$ for $t>0$. Notice that $ \bl(X^x_t) $ and $
\bl(X^x_{t^-})$ have the same distribution.\label{density-dim_inf}
\end{theorem}

\begin{proof}
The proof applies the following simple corollary of Fubini's theorem
on $\mathbb{R}^M$ with Lebesgue measure $ \lambda $ and a
probability space $(\Omega,\mathcal{F},P) $: let $ \nu $ be a
probability measure on $\mathbb{R}^M \times \Omega $ such that there
is random density $ p:\mathbb{R}^M \times \Omega \to
\mathbb{R}_{\geq 0} $ with
$$
\int_{\mathbb{R}^M \times \Omega} f(x,\eta) p(x,\eta) (\lambda\otimes P)(dx,d\eta) = \int_{\mathbb{R}^M \times \Omega} f(x,\eta) \nu(dx,d\eta),
$$
then the marginal of $ \nu $ on $\mathbb{R}^M $ is absolutely continuous with respect to Lebesgue measure $ \lambda $ with density
$$
p(x)= \int_{\Omega} p(x,\eta) P(d\eta)
$$
for almost all $ x \in \mathbb{R}^M $. In our case we know that the
law of $ X^x_t $ is absolutely continuous for almost all
trajectories $ \eta $ of the compound Poisson processes $ L$
(precisely those where $\tau_1 > 0 $ where only one coordinate jumps at
each jumping time, and finitely many jumps occur on compact intervals),
the probability measure $\nu$ corresponds to the distribution of
$(\bl(X^x_t),L) $, where we choose $\Omega$ as the space of c\`adl\`ag
trajectories on $ \mathbb{R}_{\geq 0} $ with values in $
\mathbb{R}^m $. Finally we have that the law of $ X^x_t $ is $ p(x)
\lambda(dx) $.
\end{proof}

\section{Applications of the infinite dimensional result to Interest Rate Theory}

In mathematical Finance the theory of interest rates deals with the market of interest rate related products like swaps, bonds, bills, etc. If one considers default-free products one can crystallize from the data of real markets the prices of default-free zero-coupon bonds for any maturity. A zero coupon bond contract with maturity $ T $ (a calendar date) can be entered at calendar time $ t \leq T $ and (certainly) pays $ 1 $ unit of currency at maturity time $ T $. Therefore bonds reflect the level of interest rate between time $ t $ and maturity time $ T $. No coupons are paid between $ t $ and $ T $, which explains the notion zero-coupon bond. We denote the price of a default-free zero-coupon bond with maturity $ T $ at time $ t \leq T $ by $ P(t,T) $. Commonly one assumes that bond prices are at least $ C^1 $ with respect to $ T $, i.e.
$$
P(t,T) = \exp(-\int_t^T f(t,r) dr)
$$
for $ 0 \leq t \leq T $. This leads to the concept of the short rate
$$
R_t = f(t,t)
$$
for $ t \geq 0 $, which corresponds to the level of interest rate for instantaneous transactions from $ t $ to $ t + dt $. As usual in mathematical Finance discounted default-free zero coupon bonds are modeled by semimartingales and one assumes the existence of an equivalent martingale measure for discounted price processes. This leads to the following fundamental formula with respect to the martingale measure
$$
E \bigl(\exp(-\int_t^T R_s ds)| \mathcal{F}_t \bigr) = P(t,T) = \exp(-\int_t^T f(t,r) dr).
$$
The formula simply expresses the fact that the expected valued of discounted value of the payoff $ P(T,T) = 1 $ conditional on today's information equals today's price with respect to the martingale measure. Assuming a jump diffusion model for $ {(f(t,T)}_{ 0 \leq t \leq T} $ for $ T \geq 0 $ for the forward rates together with Musiela's parametrization $ r(t,T-t) = f(t,T) $ for $ 0 \leq t \leq T $ leads to the famous Heath Jarrow Morton (HJM) equation of interest rate theory, which is an SDE taking values in a Hilbert space of forward rate curves $ H $ (and therefore an SPDE). We quote here as leading reference \cite{FilTap07}, where the no arbitrage conditions for the HJM equation are discussed in all necessary details. A very readable introduction can also be found in \cite{bjogom:99}, in particular for HJM equations with jumps.

The HJM equation has been analysed from different points of view:
\begin{itemize}
 \item The question, which HJM equations driven by finitely many Brownian motions admit finite dimensional realizations has been treated in all detail in \cite{Fil03}. This research has been inspired by \cite{BjoSve01}, where the geometric approach has been introduced. The satisfying answer is that -- under quite natural restrictions -- finite dimensional realizations do exist if and only if the corresponding factor processes are affine processes. This is the case for instance for Vasi\v{c}ek's model or for the Cox-Ingersoll-Ross model of interest rate theory. In both cases the finite dimensional realizations are in fact two dimensional. In \cite{bjogom:99} finite dimensional realizations are treated for HJM equations with jumps, however, under the strong restriction that the vector fields do not depend on the forward rate. In this case one can solve the HJM equation explicitly by variation of constants and read off the respective geometric properties of the solution process.
 \item The question whether the solution process of a HJM equation always admits a density with respect to Lebesgue measure when projected to a finite dimensional subspace has been treated in \cite{BauTei05} and could be answered affirmatively under H\"ormander type conditions.
\end{itemize}

We ask here the question -- having the theory of the previous sections in mind -- if a structure of finite dimensional realizations, such as for Vasi\v{c}ek's model, can be perturbed so strongly through the introduction of jumps that the resulting HJM evolution is ``hypo-elliptic'', i.e. the assumptions of Theorem \ref{theorem2} are fulfilled. We can answer this question affirmatively in the case of a Vasi\v{c}ek model. In contrast to \cite{bjogom:99} we allow the vector fields to be state-dependent (therefore we cannot hope for explicit solutions of the HJM equation and we have to apply local methods from differential geometry to conclude).

We consider the following HJM-model with jumps
$$
\left\{\begin{array}{ccl}
dr_t&=&\ds\left(\frac{d}{dx}r_{t^-}+\alpha_{HJM}(r_{t^-})+\beta_{HJM}(r_{t^-})\right)dt+\sigma(r_{t^-})dB_t+\delta(r_{t^-})dN_t\\[2mm]
r_0&=&\ds r^{*}\in H
\end{array}\right.
$$
on some Hilbert space $H$ of forward rate curves as constructed in \cite{BauTei05} and \cite{FilTap07}. Here $(B_t)_{t\geq 0}$ is a standard Brownian motion and $(N_t)_{t\geq 0}$ is a standard Poisson process with intensity $\widetilde{\lambda} >0  $ and jump measure $ \mu =
\delta_1 $ (hence $ \widetilde{\lambda} = \lambda $). Define
$$
\begin{array}{ccl}
\Psi_1(z)&\equiv&\ds \ln \mathbb{E}\left[e^{zW_1}\right]=\ln
\left(e^{\frac{z^2}{2}}\right)=\frac{z^2}{2}\quad\mbox{and}\\[2mm]
\Psi_2(z)&\equiv&\ds \ln \mathbb{E}\left[e^{zN_1}\right]=\ln
\left(\exp\left(\lambda\left(e^z-1\right)\right)\right)=\lambda
\left(e^z-1\right).
\end{array}
$$

Then we know from \cite{FilTap07}, equation (2.4), that
$$
\begin{array}{ccl}
\alpha_{HJM}(r)(x)&=&\ds -\sigma(r)(x)\Psi'_1\left(-\int_0^x
\sigma(r)(y)dy\right)\\[2mm]
&=&\ds \sigma(r)(x)\int_0^x \sigma(r)(y)dy
\end{array}
$$
and
$$
\begin{array}{ccl}
\beta_{HJM}(r)(x)&=&\ds -\delta(r)(x)\Psi'_2\left(-\int_0^x
\delta(r)(y)dy\right)\\[2mm]
&=&\ds -\lambda \delta(r)(x)\exp\left(-\int_0^x
\delta(r)(y)dy\right).
\end{array}
$$
For an explicit example we choose
$$
\sigma(r)(x)= \sigma(x)>0\quad\mbox{and}\quad
\delta(r)(x)=-\frac{d}{dx}\ln(B(r))(x),
$$
where the vector field $B$ will be determined later. Then we have
$$
\alpha_{HJM}(r)(x)=\sigma(x)\int_0^x \sigma(y)dy
$$
and
$$
\beta_{HJM}(r)(x)=\lambda \frac{d}{dx}\ln(B(r))(x)\frac{(B(r))(x)}{(B(r))(0)}
=\lambda \frac{\frac{d}{dx}(B(r))(x)}{(B(r))(0)}.
$$
We choose $ B $ such that $(B(r))(x)$ is positive on $H$, for $ x \in \mathbb{R} $ and $ B(r)(0) = 1 $
for all $ r \in H $, whence $\delta$ is well defined. Thus, for such
$r \in U$ we have
$$
\beta_{HJM}(r)(x) = \lambda \frac{d}{dx}(B(r))(x)
$$
and
$$
\delta(r)(x) =-\frac{d}{dx} \ln(B(r))(x).
$$
A particular choice in the spirit of Remark \ref{banach-map-vector-fields} is given through
$$
B(r)(x)= \psi(x,l(r)),
$$
where the maps $ y \mapsto \frac{d}{dx}\psi(.,y) $ and $ y \mapsto \frac{1}{\psi(.,y)} $ from $ \mathbb{R} $ to $ \dom(A^{\infty}) \subset H $ are supposed to be $ C^{\infty}$-bounded with $ \psi(0,y) = 1 $ for all $ y \in \mathbb{R} $. The map $ l $ denotes here a non-vanishing linear functional $ l : H \to \mathbb{R} $. Hence $ \delta $ and $ \beta $ are well-defined $C^{\infty}$-bounded vector fields on the whole Hilbert space and we have global existence of mild solutions.

The Vasi\v{c}ek model is defined by
$$
\sigma(r) = \rho \exp(- a x),
$$ 
for $ \rho, a > 0 $ without any jump component. By \cite{BjoSve01} and \cite{Fil03} we know that the Vasi\v{c}ek-model
admits finite dimensional realizations as for
$$
V_0(r)(x)= \frac{d}{dx}r(x)+\alpha_{HJM}(r)(x)
$$
we have
$$
\dim (\{V_{0},\sigma\}_{LA}(r))\leq 2
$$
and at any point $ r \in \dom({(\frac{d}{dx})}^{\infty}) $. Here the index $ LA $ stands for the Lie algebra generated by the vector fields $ V_0,V_1 $ on $ \dom(A^{\infty}) $. If we add a jump structure as described above and if we choose $ \psi $ generic, the two dimensional structure (a regular finite dimensional realization in the sense of \cite{Fil03}) is destroyed, since then the drift changes due to no-arbitrage. We obtain a dense Lie algebra if we choose the vector field $ B $ generically.

These results might be of interest for recent works in interest rate theory, see for instance \cite{malmanrec:07}, where under diffusion assumptions hypo-ellipticity has been tested empirically. If one allows for jumps in an HJM model the phenomenon of hypo-ellipticity seems to be more generic.

\section{Smooth densities for the law $X_{T}^{x}$ on $ \mathbb{R}^M$}

In the sequel we consider the case $\dim H = M $ and $ \bl =
\operatorname{id} $ and we choose a coordinate representation $ H =
\mathbb{R}^M $. We then want to show that the $p$-th power of the
inverse of the Malliavin covariance matrix of $ X^x_t $ for $ t >0 $
can be integrated even with respect to Poisson trajectory $ \eta $.
We therefore need an extension of the H\"ormander condition, which
is called the uniform H\"ormander condition:

Following \cite{Nua95}, we define
\begin{align*}
\Sigma_{0}^{\prime}  &  :=\{V_{1},\dots,V_{d}\}\\
\Sigma_{n}^{\prime}  &  :=\left\{  [V_{k},V],k=1,\dots,d,V\in\Sigma
_{n-1}^{\prime};[V_0,V]+\frac{1}{2}\sum_{i=1}^{d}[V_{i,}[V_{i}
,V]],V\in\Sigma_{n-1}^{\prime}\right\}
\end{align*}
for $n\geq1$. We assume that there exists
$j_{0}$ and $ c>0$ such that
\begin{equation}
\inf_{\xi\in S^{M-1}}\sum_{j=0}^{j_{0}}\sum_{V\in\Sigma_{j}^{\prime}
}\left\langle V(x),\xi\right\rangle ^{2}\geq c \label{uniformhoermander}
\end{equation}
uniformly in $x\in\mathbb{R}^{M}$.

\begin{theorem}
Assume that $ \dim H < \infty $. We take Assumptions \ref{ass
hoermander condition} and \ref{ass jump condition} for granted, but
assume that the H\"ormander condition (\ref{hoermander condition})
holds true uniformly on $\mathbb{R}^M$ in the sense of \eqref{uniformhoermander}. Let ${(X_{t}^{x})}_{t\geq0}$
denote the unique c\`{a}dl\`{a}g solution of equation (\ref{equ1})
and fix $ t >0$. Then the random variable $X_{t}^{x}$ admits a
smooth density with respect to Lebesgue measure on $\mathbb{R}^M$.
Furthermore, the covariance matrix of $ X^x_t$ is invertible with
$p$-integrable inverse for all $p\geq1$.\label{smooth density}
\end{theorem}

\begin{proof}
We write the Malliavin derivative of $X_{t}^{x}$,
\[
D_{s}^{i}X_{t}^{x}=J_{0\to t}(x)J_{0\to s}(x)^{-1}
V_{i}(X_{s^-}^{x})\mathbf{1}_{[0,t]}(s).
\]
and calculate the reduced covariance matrix
\[
\left\langle C_t\xi,\xi\right\rangle
=\sum_{i=1}^{d}\int_{0}^{t}\left\langle J_{0\to
s}(x)^{-1}V_{i}(X_{s^-}^{x}),\xi\right\rangle ^{2}ds.
\]
We now apply the result from Theorem \ref{structure theorem} and condition on
the trajectories of the compound Poisson process (\ref{poisson-process})
\begin{multline}
\sup_{\xi\in S^{M-1}}P(\left\langle C_{t}\xi,\xi\right\rangle <\epsilon)\\
=\sup_{\xi\in S^{M-1}}\sum_{n_{1},\dots,n_{m} \geq
0}\Big[\prod_{k=1}^{m}P(N_{t} ^{k}=n_{k})\Big]P(\left\langle
C_{t}\xi,\xi\right\rangle <\epsilon|N_{t}^{j} =n_{j}\text{ for
}j=1,\dots,m).\label{decomposition}
\end{multline}
As in the proof of Theorem \ref{theorem2}, we can decompose
$\langle C_t \xi, \xi\rangle$ into
\[
\left\langle C_{t}\xi,\xi\right\rangle =\sum_{k=0}^{\infty}\sum_{i=1}^{d}
\int_{\tau_{k}\land t}^{\tau_{k+1}\land t}\left\langle J_{0\to s}(x)^{-1}
V_{i}(X_{s^-}^x),\xi\right\rangle ^{2}ds,
\]
where $\tau_{0}=0<\tau_{1}< \ldots <\tau_{n}\leq\dots$ denotes the
sequence of jump times of $(N_{t})_{0\leq t\leq T}$. Hence, we
obtain for $n=n_{1}+\dots+n_{m}$
\begin{multline*}
\sup_{\xi\in S^{M-1}}P(\left\langle C_{t}\xi,\xi\right\rangle
<\epsilon
|N_{t}^{j}=n_{j},\ j=1,\dots,m) \leq \\
\sup_{\xi\in S^{M-1}}P\left(   \sum_{i=1} ^{d}\int_{\tau_{k}\land
t}^{\tau_{k+1}\land t}\left\langle J_{0\to
s}(x)^{-1}V_{i}(X_{s^-}^{x}),\xi\right\rangle ^{2}ds<\epsilon \Big
\vert N_{t}^{j}=n_{j},\ j=1,\dots,m\right),
\end{multline*}
for all $ 0 \leq k \leq n $. Observing that $\max_{0\leq k\leq n}
(\tau_{k+1}-\tau_{k})^{K(p)}\geq(\frac{t}{n})^{K(p)}$ (after all
we only have $n$ jumps, so the maximal distance between two
consecutive jumps is bigger than $\frac{t}{n}$), we finally obtain
\[
P(\left\langle C_{t}\xi,\xi\right\rangle
<\epsilon|N_{t}^{j}=n_{j}\text{ for }j=1,\dots,m)\leq\epsilon^{p}
\]
for $0\leq\epsilon\leq(\frac{t}{n})^{K(p)}\epsilon_{0}(p)$ due to
the calculations outlined in the Appendix. Note that we can apply the
calculations from the Appendix, since $J_{0\to s}(x)^{-1}$
is well-defined and bounded due to boundedness of
$(id+zd\delta_{j})^{-1}$ for $ z \in \operatorname*{supp}(\mu_j) $
and $j=1,\dots,m$. Hence integration with respect to the measures
$\mu_{j}$ is possible and yields finite bounds. Recall also that
$\mu_{j}$ has moments of all orders, hence $X_{t}^{x}$ is $L^{p}$
and so is $J_{0\to s}(x)^{-1}$ (see \cite{Pro90} for all necessary
details on SDEs).

Let $\Lambda = \inf_{\xi\in S^{M-1}}\left\langle C_{t}\xi,\xi\right\rangle$ be
the smallest eigenvalue of the reduced covariance matrix $C_t$. Following the
steps of \cite{Nua95}, Lemma 2.3.1, we know that
\[
P(\Lambda <\epsilon \; |\; N_{t}^{j}=n_{j}\text{ for
}j=1,\dots,m)\leq \rm{const}\cdot\epsilon^{p}
\]
for any $p\geq2$ and $0\leq\epsilon\leq
\big(\frac{t}{n}\big)^{K(p+2M)}\epsilon _{0}(p+2M) =:
\epsilon_{max}$, where the constant depends on the $p$-norm of
$C_{t}$. In the sequel we shall denote any constant of this type by
$ D $. We denote by $ \rho $ the law of $ \Lambda $ conditioned on $
N^j_t = n_j $ for $ j=1,\ldots,m $. Consequently, for $j = 1,
\ldots, m$, we have by Fubini's Theorem
\begin{equation*}
\begin{split}
E\left(\frac{1}{\Lambda ^{p-1}}\; \left|\;
N_{t}^{j}=n_{j}\right.\right) &= E\left(\frac{1}{\Lambda
^{p-1}}\cdot \mathbf{1}_{\{\Lambda > \epsilon_{max}\}}\; \left|\;
N_{t}^{j}=n_{j}\right.\right) \\ & \quad + E\left(\frac{1}{\Lambda
^{p-1}}\cdot \mathbf{1}_{\{\Lambda \leq \epsilon_{max}\}}\; \left|\;
N_{t}^{j}=n_{j}\right.\right)\\ &\leq \frac{1}{\epsilon_{max}^{p-1}}
+ \int_0^{\epsilon_{max}} \frac{1}{z^{p-1}} \rho(dz)\\ & =
\frac{1}{\epsilon_{max}^{p-1}} + \int_0^{\epsilon_{max}} (p-1)
\int_z^{\infty} \frac{1}{t^p} dt \rho(dz)
\\ &=\frac{1}{\epsilon_{max}^{p-1}} + (p-1)\int_0^{\epsilon_{max}}
\frac{1}{z^{p}} \int_{0}^{z} \rho(dt) dz + \\ & +
(p-1)\int_{\epsilon_{max}}^{\infty} \frac{1}{z^{p}} \int_{0}^{\epsilon_{max}}
\rho(dt) dz \\ & \leq \frac{D}{\epsilon_{max}^{p-1}} +
D\underbrace{\int_0^{\epsilon_{max}} \frac{1}{z^{p}}\; z^p dz}_{=
\epsilon_{max}}\\ & \leq D (\frac{t}{n})^{K(p+2M)}\epsilon_{0}(p+2M)+
\\ &+ \frac{D}{[(\frac{t}{n})^{K(p+2M)}\epsilon_{0}(p+2M)]^{p-1}}.
\end{split}
\end{equation*}
Here we applied $ \int_0^{z} \rho(dt) \leq \rm{const} z^p $ as previously proved. Hence through the decomposition (\ref{decomposition}),
\begin{multline*}
E\left(\frac{1}{\Lambda^{p-1}
}\right)\leq\sum_{n_{1},\dots,n_{m}>0}\prod_{k=1}^{m}P(N_{t}^{k}=n_{k})\cdot\\
\cdot D \cdot
\left[\Big(\frac{t}{n}\Big)^{K(p+2M)}\epsilon_{0}(p+2M)+\frac{1}{[(\frac{t}
{n})^{K(p+2M)}\epsilon_{0}(p+2M)]^{p-1}}\right]  <\infty
\end{multline*}
the result follows by $n=n_{1}+\dots+n_{m}$ and by the following
fact for any real number $ K $,
\[
\sum_{n_{1},\dots,n_{m}>0}\frac{\widetilde{\lambda}_{1}^{n_{1}}\cdots\widetilde{\lambda}_{m}^{n_{m}}}{n_{1}!\cdots n_{m}!}e^{-t \widetilde{\lambda}_{1}n_{1}- \dots - t \widetilde{\lambda}_{m} n_{m}}t^n (n_{1}+\dots + n_{m})^{K} < \infty. \qedhere\]
\end{proof}
\begin{remark}
We could have also applied the beautiful results of \cite[Corollary 3.25]{kusstr:85} to evaluate the $ L^p$ norm of the inverse of the covariance matrix between two jumps. Both ways lead to the same result. We have been choosing our approach since we could it root it is much as possible in the standard reference \cite{Nua95}.
\end{remark}

\section{Calculating the Greeks in finite dimension}

In the sequel we consider the case $\dim H = M $ and $ \bl =
\operatorname{id} $ as in the previous section. Once we are given an
invertible Malliavin covariance matrix with $p$-integrable inverse
such as in Theorem \ref{smooth density}, we can easily calculate
derivatives with respect to initial values and obtain explicit
formulas for so-called Malliavin weights (see \cite{Fou99} for successful applications of this method in mathematical Finance). We sum up quickly the main idea: in mathematical Finance the gradient of the function $ x \mapsto E(f(X^x_t)) $ has the meaning of hedging ratios, which control the hedging portfolios away from jumps. Hence for any hedging portfolio corresponding to prices $ E(f(X^x_t)) $ of a certain derivative at maturity $ t > 0 $ it is crucial to know $ \nabla E(f(X^x_t)) $ to perform hedging off jumps.

Very often pricing results in the applications of weak-approximation-scheme for the process $ X $, for instance the Euler-Maruyama scheme. For the calculation of $ \nabla E (f(X^x_t)) $ in the direction of some vector $ v \in H $ basically three methods can be applied:
\begin{itemize}
 \item a finite difference method to approximate $ \nabla E (f(X^x_t)) \cdot v  $ resulting in the calculation of $ \frac{E(f(X^{x+\epsilon v}_t)) - E(f(X^{x}_t))}{\epsilon} $ ($v \in H$ denotes some vector) for small $ \epsilon > 0 $.
 \item a pathwise method applying the formula
$$
\nabla E (f(X^x_t)) \cdot v = E( d f(X^x_t) J_{0 \to t}(x) \cdot v)
$$
resulting in the weak numerical approximation of $ (X^x_t,J_{0 \to t}(x)) $.
\item the method of Malliavin weights applying the formula
$$
\nabla E (f(X^x_t)) \cdot v = E( f(X^x_t) \pi^{v})
$$
results in the weak numerical approximation of $ (X^x_t,\pi^{v}) $.
\end{itemize}

The first method is the most robust one in the sense that it can be applied under very weak assumptions both on $ X^x_t $ and on the payoff $ f$, but the rate of convergence might be very slow since the errors of Monte-Carlo evaluations are amplified. The second method works for all reasonable jump-diffusion processes but one needs Lipschitz conditions on the payoff $ f $. The third method needs the assumptions of Theorem \ref{smooth density} on $ X^x_t $, but no restrictions on the payoff $ f $, which makes the third method attractive for several problems from mathematical Finance, where measurable, non-Lipschitz payoffs (e.g.~digital options) are quite usual and hypo-ellipticity assumptions as in Theorem \ref{smooth density} are common, too.

The implementation of procedures for all three methods have been outlined in \cite{gobmun:05} in the pure diffusion case. We shall not work on this issue here, since our main message is that one can implement precisely the same methods as in pure diffusion cases for jump-diffusions. The important point is that the formulas have the same structure in both cases, a fact, on which we shall point in this section at several occasions.

We denote in this sequel the Skorohod integral (resp. the divergence operator) by $\delta$ and its domain by $\dom(\delta).$

\begin{definition}
 Assume that $ H=\mathbb{R}^M$, fix $t>0$ and a direction $v\in\mathbb{R}^{M}$. We define a set of Skorohod-integrable processes
\[
\mathbb{A}_{t,x,v}=\bigg\{a\in\operatorname*{dom}(\delta)\text{ such that }
\sum_{i=1}^{d}\int_{0}^{t}J_{0\to s}(x)^{-1}V_{i}(X_{s^-}^{x})a_{s}
^{i}ds=v\bigg\}
\]
and call it the set of path-perturbations with target-value $ v $.
\end{definition}
\begin{remark}
In the previous definition such as in the whole section assertions on Skorohod-integrability are meant Poissonian-trajectory-wise.
\end{remark}
\begin{proposition}
Assume that $H =\mathbb{R}^M$. We take Assumption \ref{ass jump
condition} for granted. Fix $t>0$ and a direction
$v\in\mathbb{R}^{M}$. Assume furthermore uniform ellipticity, i.e.
$M=d$ and there is $c>0$ such that
\[
\inf_{\xi\in S^{M-1}}\sum_{k=1}^{M}\left\langle V_{k}(x),\xi\right\rangle
^{2}\geq c.
\]
Then $\mathbb{A}_{t,x,v}\neq\emptyset$ and there exists an
integrable, real valued random variable $\pi^{v}$ (which depends
linearly on $v$) such that for all bounded
random variables $f$ we obtain
\[
\frac{d}{d\epsilon} \vert _{\epsilon=0}E(f(X_{t}^{x+\epsilon
v}))=E(f(X_{t}^{x})\pi^{v}).
\]
Such a random variable $ \pi^{v} $ is called a Malliavin weight and can be obtained through an It\^o integral.
\end{proposition}

\begin{remark}
The assertion of this theorem corresponds to Assumption (E) in \cite{gobmun:05} and to the assumptions of \cite{Fou99}. The assumptions are seen as too restrictive since not every problem in mathematical Finance has an elliptic volatility matrix. The formulas of \cite{gobmun:05} and \cite{Fou99} correspond precisely to the formulas obtained here, which leads to the assertion that even in the presence of jumps one can apply the same (numerical) methods for the calculation of greeks.
\end{remark}

\begin{proof}
Here the proof is particularly simple, since we can take a matrix
$\sigma(x):=(V_{1}(x),\dots,V_{M}(x))$, which is uniformly invertible with
bounded inverse. We define
\[
a_{s}:=\frac{1}{t}\sigma(X_{s^-}^{x})^{-1}\cdot {J_{0 \to s}(x)} \cdot v
\]
for $0\leq s\leq t$ and obtain that $a\in\mathbb{A}_{t,x,v}$. Furthermore --
as in \cite{Fou99} and \cite{DavJoh06} -- we obtain
\[
\pi^{v}=\sum_{i=1}^{M}\int_{0}^{t}a_{s}^{i}dB_{s}^{i},
\]
since the Skorohod integrable process $a$ is in fact adapted, left-continuous
and hence It\^{o}-integrable.
\end{proof}

\begin{theorem}
Assume that $H =\mathbb{R}^M$. We take Assumptions \ref{ass
hoermander condition} and \ref{ass jump condition} for granted, but
assume that the H\"ormander condition (\ref{hoermander condition})
holds true uniformly on $\mathbb{R}^M$ (see Section 7). Fix $t>0$ and a direction
$v\in\mathbb{R}^{M}$. Then $\mathbb{A}_{t,x,v}\neq\emptyset$ and
there exists an integrable, real valued random variable $\pi^{v}$ (which
depends linearly on $v$) such that for all bounded random variables
$f$ we obtain
\[
\frac{d}{d\epsilon} \vert _{\epsilon=0}E(f(X_{t}^{x+\epsilon
v}))=E(f(X_{t}^{x})\pi^{v}).
\]
We can choose $\pi^{v}$ to be the Skorohod integral of any element $a\in
\mathbb{A}_{t,x,v}\neq\emptyset$ and call it a Malliavin weight. Moreover, by the explicit construction of $ a $ in the proof, we can assert that $ \pi^v $ is the sum of an It\^o integral and an integral with respect to Lebesgue measure, see for instance \cite{gobmun:05}.
\end{theorem}

\begin{remark}
The assertion of this theorem corresponds to Assumption (E') in \cite{gobmun:05}. The assumptions (E) and (E') are fundamental for the third method in \cite{gobmun:05}. Again the formulas of \cite{gobmun:05} correspond to the formulas obtained here.
\end{remark}

\begin{proof}
We take $f$ bounded with bounded first derivative, then we obtain
\[
\frac{d}{d\epsilon} \vert _{\epsilon=0}E(f(X_{t}^{x+\epsilon
v}))=E(df(X_{t}^{x})J_{0\to t}(x)\cdot v).
\]
If there is $a\in\mathbb{A}_{t,x,v}$, we obtain
\begin{align*}
E(df(X_{t}^{x})J_{0\to t}(x)\cdot v)  &  =E(df(X_{t}^{x})\sum
_{i=1}^{d}\int_{0}^{t}J_{0\to t}(x)J_{0\to s}(x)^{-1}
V_{i}(X_{s^-}^{x})a_{s}^{i}ds)\\
&  =E(\sum_{i=1}^{d}\int_{0}^{t}df(X_{t}^{x})J_{0\to t}
(x)J_{0\to s}(x)^{-1}V_{i}(X_{s^-}^{x})a_{s}^{i}ds)\\
&  =E(\sum_{i=1}^{d}\int_{0}^{t}D_{s}^{i}f(X_{t}^{x})a_{s}^{i}ds)\\
&  =E(f(X_{t}^{x})\delta(a)).
\end{align*}
Here we cannot assert that the strategy is It\^{o}-integrable, since it will be
anticipative in general. In order to see that $\mathbb{A}_{t,x,v}\neq
\emptyset$ we construct an element, namely
\[
a_{s}^{i}:=\left\langle J_{0\to s}(x)^{-1}V_{i}(X_{s^-}^{x}
),(C_{t})^{-1}v\right\rangle ,
\]
where $C_{t}$ denotes the reduced covariance matrix from Theorem \ref{smooth
density}. Indeed
\begin{align*}
&  \sum_{i=1}^{d}\left\langle \int_{0}^{t}J_{0\to s}(x)^{-1}
V_{i}(X_{s^-}^{x})a_{s}^{i}ds,\xi\right\rangle \\
&  =\sum_{i=1}^{d}\int_{0}^{t}\left\langle J_{0\to s}(x)^{-1}
V_{i}(X_{s^-}^{x}),\xi\right\rangle \left\langle J_{0\to s}
(x)^{-1}V_{i}(X_{s^-}^{x}),(C_{t})^{-1}v\right\rangle ds\\
&  =\left\langle \xi,C_{t}(C_{t})^{-1}v\right\rangle =\left\langle
\xi,v\right\rangle
\end{align*}
for all $\xi\in\mathbb{R}^{M}$, since $C_{t}$ is a symmetric
random operator
defined via
\[
\left\langle \xi,C_{t}\xi\right\rangle
=\sum_{i=1}^{d}\int_{0}^{t}\left\langle J_{0\to
s}(x)^{-1}V_{i}(X_{s^-}^{x}),\xi\right\rangle ^{2}ds
\]
for $\xi\in\mathbb{R}^{M}$.
\end{proof}

For any other derivative with respect to parameters $\epsilon$, we
consider a modified set, namely
\[
\mathbb{B}_{t,x,v}=\bigg\{b\in\operatorname*{dom}(\delta)  \vert \, \sum
_{i=1}^{d}\int_{0}^{t}J_{0\to s}(x)^{-1}V_{i}(X_{s^-}^{x})b_{s}
^{i}ds={J_{0\to t}(x)}^{-1} \frac{d}{d\epsilon}\vert
_{\epsilon=0}X_{t}^{x,\epsilon}\bigg\}.
\]
Here we are given a parameter-dependent process
$X_{t}^{x,\epsilon}$, where all derivatives with respect to
$\epsilon$ can be calculated nicely. Also in this case we can
construct -- if the reduced covariance matrix is invertible and regular enough -- an
element, namely
\[
b_{s}^{i}:=\left\langle J_{0\to
s}(x)^{-1}V_{i}(X_{s^-}^{x}),(C_{t})^{-1}{J_{0\to t}(x)}^{-1} \frac{d}{d\epsilon} \vert_{\epsilon=0}X_{t}^{x,\epsilon}\right\rangle .
\]
This is a consequence of the following reasoning,
\begin{align*}
& \sum_{i=1}^{d}\Big\langle \int_{0}^{t}J_{0\to s}(x)^{-1}V_{i}
(X_{s^-}^{x})b_{s}^{i}ds,\xi\Big\rangle =\sum_{i=1}^{d}\int_{0}^{t}\Big\langle J_{0\to s}(x)^{-1}V_{i}
(X_{s^-}^{x}),\xi\Big\rangle \cdot \\
& \hspace{2mm} \cdot\Big \langle J_{0\to s}(x)^{-1}
V_{i}(X_{s^-}^{x}),(C_{t})^{-1}{J_{0\to t}(x)}^{-1} \frac{d}{d\epsilon }\vert_{\epsilon=0}X_{t}^{x,\epsilon}\Big\rangle ds\\ & =\Big\langle
\xi,C_{t}(C_{t})^{-1}{J_{0\to t}(x)}^{-1} \frac
{d}{d\epsilon} \vert_{\epsilon=0}X_{t}^{x,\epsilon}\Big\rangle =\Big\langle
\xi,{J_{0\to t}(x)}^{-1} \frac{d}{d\epsilon}\vert_{\epsilon=0}X_{t}^{x,\epsilon}\Big\rangle ,
\end{align*}
due to symmetry of $C_{t}$.

\section{Appendix}

\begin{theorem}
Let $(\Omega,\mathcal{F},P,(\mathcal{F}_{t})_{t\geq0})$ be a
filtered probability space and let $(B_{t})_{t\geq0}$ be a
$d$-dimensional Brownian motion adapted to the filtration (which
is not necessarily generated by the Brownian motion). Let
$V,V_{1},\dots,V_{d}$, the diffusion vector fields be
$C^{\infty}$-bounded on $\mathbb{R}^{M}$ and consider the
continuous solution $(X_{t}^{x})_{0\leq t\leq T}$ of a stochastic
differential equation (in
Stratonovich notation). $V_{0}$ denotes the Stratonovich corrected drift term,
\begin{align}
dX_{t}^{x}  &  =V_{0}(X_{t}^{x})dt+\sum_{i=1}^{d}V_{i}(X_{t}^{x})\circ
dB_{t}^{i},\\
X_{0}^{x}  &  =x\text{.}
\end{align}
Assume that the uniform H\"ormander condition holds true (see the
proof for the precise statement). Then for any $p\geq1$ there exist
numbers $\epsilon_{0}(p)>0$ and an integer $K(p)\geq1$ such that for
each $0<t<T$
\[
\sup_{\xi\in S^{M-1}}P(\left\langle C_{t}\xi,\xi\right\rangle
<\epsilon
)\leq\epsilon^{p}
\]
holds true for $0\leq\epsilon\leq t^{K(p)}\epsilon_{0}(p)$. The
result holds uniformly in $x$.
\end{theorem}
\begin{remark}
The time-dependence of the estimate $0\leq\epsilon\leq t^{K(p)}\epsilon_{0}(p)$ is best explained by re-doing the proof. It is heavily applied in Section 7 and the main technical ingredient of the given proof. We could have also used directly the results from \cite{kusstr:85}.
\end{remark}
\begin{proof}
The proof of the theorem is a careful re-reading of the Norris Lemma and the
classical proof of the H\"ormander theorem in probability theory (see
\cite{Mal97} or \cite{Nua95}). We shall sketch this path in the sequel
(see \cite{Nua95}, pp.120--123):

\begin{enumerate}
\item Consider the random quadratic form
\[
\left\langle C_{t}\xi,\xi\right\rangle
=\sum_{i=1}^{d}\int_{0}^{t}\left\langle J_{0\to
s}(x)^{-1}V_{i}(X_{s}^{x}),\xi\right\rangle ^{2}ds.
\]
Following \cite{Nua95}, we define
\begin{align*}
\Sigma_{0}^{\prime}  &  :=\{V_{1},\dots,V_{d}\}\\
\Sigma_{n}^{\prime}  &  :=\left\{  [V_{k},V],k=1,\dots,d,V\in\Sigma
_{n-1}^{\prime};[V_{0},V]+\frac{1}{2}\sum_{i=1}^{d}[V_{i,}[V_{i}
,V]],V\in\Sigma_{n-1}^{\prime}\right\}
\end{align*}
for $n\geq1$. We assume that there exists
$j_{0}$ and $ c>0$ such that
\[
\inf_{\xi\in S^{M-1}}\sum_{j=0}^{j_{0}}\sum_{V\in\Sigma_{j}^{\prime}
}\left\langle V(x),\xi\right\rangle ^{2}\geq c
\]
uniformly in $x\in\mathbb{R}^{M}$.

\item We define $m(j):=2^{-4j}$ for $0\leq j\leq j_{0}$ and the sets
\[
E_{j}:=\Big\{\sum_{V\in\Sigma_{j}^{\prime}}\int_{0}^{t}\left\langle
J_{0\to s}(x)^{-1}V(X_{s}^{x}),\xi\right\rangle ^{2}ds\leq
\epsilon^{m(j)}\Big\}.
\]
We consider the decomposition
\begin{align*}
E_{0}  &  =\{\left\langle C_{t}\xi,\xi\right\rangle \leq\epsilon
\}\subset(E_{0}\cap E_{1}^{c})\cup(E_{1}\cap
E_{2}^{c})\cup\dots\cup
(E_{j_{0}-1}\cap E_{j_{0}}^{c})\cup F,\\
F  &  =E_{0}\cap\dots\cap E_{j_{0}}.
\end{align*}
and proceed with
\[
P(F)\leq C\epsilon^{\frac{q\beta}{2}},
\]
for $\epsilon\leq\epsilon_{1}$ and any $q\geq 2$ with a constant $C$
depending on $q$ and the norms of the derivatives of the vector
fields $V_{0},\dots,V_{d}$. Furthermore $0<\beta<m(j_{0}).$ The
number $\epsilon_{1}$ is determined by the following two (!)
equations
\begin{align*}
(j_{0}+1)\epsilon_{1}^{m(j_{0})}  &  <\frac{c\epsilon_{1}^{\beta}}{4},\\
\epsilon_{1}^{\beta}  &  <t.
\end{align*}
Hence $\epsilon_{1}$ depends on $j_{0}$, $c$, $t$ and the choice of $\beta$,
via
\[
\epsilon_{1}<\min\left(  t^{\frac{1}{\beta}},\left(  \frac{c}{4(j_{0}
+1)}\right)  ^{\frac{1}{m(j_{0})-\beta}}\right).
\]
This little observation additional to the proof in \cite {Nua95} is key for
our proof.
\item We obtain furthermore that
\begin{multline*}
P(E_{j}\cap E_{j+1}^{c}) = P \bigg( \sum_{V \in \Sigma_j^{\prime}}
\int_0^t \langle J_{0 \to s}(x)^{-1} V(X_s^x), \xi
\rangle^2 d s \leq
\epsilon^{m(j)}, \\
\sum_{V \in \Sigma_{j+1}^{\prime}} \int_0^t \langle J_{0
\to s}(x)^{-1} V(X_s^x), \xi \rangle^2 d s >
\epsilon^{m(j+1)} \bigg)\\
\leq\sum_{V\in\Sigma_{j}^{\prime}}P \Bigg(
\int_{0}^{t}\left\langle J_{0\to
s}(x)^{-1}V(X_{s}^{x}),\xi\right\rangle ^{2}ds\leq
\epsilon^{m(j)}, \\
 \sum_{k=1}^{d}\int_{0}^{t}\left\langle J_{0\to s}
(x)^{-1}[V_{k},V](X_{s}^{x}),\xi\right\rangle ^{2}ds+\\
+\int_{0}^{t}\Big\langle J_{0\to s}(x)^{-1}\Big(
[V_{0},V]+\frac{1}{2}\sum_{i=1}^{d}[V_{i},[V_{i},V]]\Big)  (X_{s}^{x}
),\xi\Big \rangle ^{2}ds>\frac{\epsilon^{m(j+1)}}{n(j)}\Bigg)  ,
\end{multline*}
where $n(j)=\#\Sigma_{j}^{\prime}$. Since we can find the bounded
variation and the quadratic variation part of the martingale
$(\left\langle J_{0\to s}(x)^{-1}V(X_{s}^{x}),\xi\right\rangle
)_{0\leq s\leq t}$ in the above expression, we are able to apply
Norris Lemma (see \cite{Nua95}, Lemma 2.3.2). We observe that
$8m(j+1)<m(j)$, hence we can apply it with $q=\frac{m(j)}{m(j+1)}$.

\item We obtain for $p\geq2$ -- still by the Norris Lemma -- the estimate
\[
P(E_{j}\cap E_{j+1}^{c})\leq d_{1}\left(  \frac{\epsilon^{m(j+1)}}
{n(j)}\right)  ^{rp}+d_{2}\exp\left(
-\Big(\frac{\epsilon^{m(j+1)}}{n(j)}\Big)^{-\nu }\right)
\]
for $\epsilon\leq\epsilon_{2}$. Furthermore $r,\nu>0$ with $18r+9\nu<q-8$, the
numbers $d_{1},d_{2}$ depend on the vector fields $V_{0},\dots,V_{d}$, and on
$p$, $T$. The number $\epsilon_{2}$ can be chosen as $\epsilon_{2}
=\epsilon_{3}t^{k_{1}}$, where $\epsilon_{3}$ does not depend on
$t$ anymore.

\item Putting all together we take the minimum of $\epsilon_{1}$
and $\epsilon_{2}$ to obtain the desired dependence on $t$.
\qedhere
\end{enumerate}
\end{proof}


\end{document}